\documentclass[review]{elsarticle}

\usepackage{lineno,hyperref}
\usepackage{amssymb,enumerate}
\usepackage{amsmath}
\usepackage{graphics}
\def\R{{\bf R}}

\def\e{{\varepsilon}}

\newtheorem{thm}{Theorem}[section]
 
 \newtheorem{lem}[thm]{Lemma}

 \newtheorem{rem}[thm]{Remark}

\journal{Journal of \LaTeX\ Templates}

\begin{document}

\begin{frontmatter}

\title{The sharp lifespan estimate for semilinear damped wave equation with Fujita critical power in higher dimensions}

\author[mymainaddress,secondaryaddress]{Ning-An Lai\corref{mycorrespondingauthor}}
\cortext[mycorrespondingauthor]{Corresponding Author}
\ead{hyayue@gmail.com}

\author[secondaryaddress]{Yi Zhou}
\ead{yizhou@fudan.edu.cn}

\address[mymainaddress]{Institute of Nonlinear Analysis and Department of Mathematics,\\ Lishui University, Lishui 323000, China}
\address[secondaryaddress]{School of Mathematical Sciences, Fudan University, Shanghai 200433, China}

\begin{abstract}

This paper is concerned with the lifespan estimate of classical solutions with small initial data to the Cauchy problem of semilinear damped wave equations with
the Fujita critical exponent. We establish the following sharp upper bound of the lifespan
\begin{equation}\nonumber\\
\begin{aligned}
T(\varepsilon)\leq \exp\left(C\varepsilon^{-\frac 2n}\right)
\end{aligned}
\end{equation}
in higher dimensions$(n\geq 4)$, by using the heat kernel as the test function. Then, together with the previous results, a complete result on the sharp lower and upper bound estimates is obtained in this case.

 \noindent {\bf R\'esum\'e}
  \newline
  \indent Cet article porte sur l'estimation de la dur\'{e}e de vie des solutions classiques avec petites donn\'{e}es initiales pour le probl\`{e}me de Cauchy de l'\'{e}uation des ondes dissipative semi-lin\'{e}aire avec l'exposant critique de Fujita. Nous \'{e}tablissons la borne sup\'{e}rieure de la dur\'{e}e de vie suivante
\begin{equation}\nonumber\\
\begin{aligned}
T(\varepsilon)\leq \exp\left(C\varepsilon^{-\frac 2n}\right)
\end{aligned}
\end{equation}
en dimension sup\'{e}rieure $(n\geq 4)$, en utilisant le noyau de la chaleur comme fonction de test. Donc, avec les r\'{e}sultats pr\'{e}c\'{e}dents, nous obtenons un r\'{e}sultat complet sur les estimations optimale de borne sup\'{e}rieure et inf\'{e}rieure
de la dur\'{e}e de vie dans ce cas.

\end{abstract}


\end{frontmatter}

\section{Introduction}

We consider the Cauchy problem of semilinear damped wave equations with Fujita critical power in $n$ dimensions:
\begin{equation}\label{1}
\left \{
\begin{aligned}
&u_{tt}-\Delta u+u_t=|u|^p,~~~(t, x)\in [0, T(\varepsilon))\times \mathbb{R}^n,\\
&u(0, x)=\varepsilon f(x),~~ u_t(0, x)=\varepsilon g(x),~~x\in \mathbb{R}^n,\\
\end{aligned} \right.
\end{equation}
where $p=1+\frac2n$, $n\geq 4$ and $\varepsilon>0$ is a parameter which represents the smallness of the initial data. We are devoted to establish the lifespan estimate from above in the form $T(\varepsilon)\leq \exp(C\varepsilon^{-\frac 2n})$, where $C$ is a positive constant independent of $\varepsilon$, then by combining the result about the lifespan estimate from below obtained in Li \cite{Li1} and Ikeda and Ogawa \cite{Ikeda}, we get the sharpness of the lifespan estimate to Cauchy problem \eqref{1}.

The study on Cauchy problem \eqref{1} with small initial data has a long history. In 1995, Li and Zhou \cite{Li} studied Cauchy problem \eqref{1} with $p$
satisfying $1<p\leq 1+\frac 2n$ in lower dimensional cases, i.e., $n=1, 2$, and established the sharp upper bound of the lifespan in the form
\begin{equation}\label{2}
T(\varepsilon)\leq\left \{
\begin{aligned}
&\exp(C\varepsilon^{-\frac2n}),~~~p=1+\frac2n,\\
&C\varepsilon^{-\frac{1}{\frac{1}{p-1}-\frac n2}},~~~~~1<p<1+\frac2n,\\
\end{aligned} \right.
\end{equation}
where $C=C(n, p, f, g)$ is a positive constant independent of $\varepsilon$. One year later Li \cite{Li1} considered the lower bound of lifespan estimate for Cauchy problem \eqref{1} for integer $p$ and $n\geq 1$ and obtained the following result by using the global iteration method introduced by Li
and Yu \cite{Li2}(see also Li and Zhou \cite{Li3}):
\begin{equation}\label{3}
T(\varepsilon)\geq\left \{
\begin{aligned}
&+\infty,~~~p>1+\frac2n,\\
&\exp(C\varepsilon^{-\frac2n}),~~~~p=1+\frac2n,\\
&C\varepsilon^{-\frac{1}{\frac{1}{p-1}-\frac n2}},~~~~~1<p<1+\frac2n,\\
\end{aligned} \right.
\end{equation}
where $C=C(n, p, f, g)$ is a positive constant independent of $\varepsilon$.
Then Todorova and Yordanov \cite {Todorova} found that Cauchy problem \eqref{1} with small initial data admits a critical power $p=1+\frac 2n$ such that Cauchy problem \eqref{1} has global solutions if $p>1+\frac2n$,
while the solutions blow up in a finite time if $1<p<1+\frac2n$. It is interesting to see that this critical power is exactly the same as the Fujita critical power for the
corresponding semilinear heat equations $v_t-\Delta v=|v|^p$, see \cite{Fujita} for details. Later, Zhang \cite{Zhang} showed that the solutions also blow up in a finite time when $p=1+\frac 2n$, by using the test function method. Nishihara \cite{Nishihara} studied the 3-D case
and established the sharp upper bound of the lifespan estimate. Up to now, the sharp lifespan
estimate from above for $p=1+\frac2n$ in higher dimensional spaces$(n\geq 4)$ is still open. Both the methods used in \cite{Li} and \cite{Nishihara} do not work in this case, since their proofs are based on the explicit formula of the solutions and the positivity of the fundamental solution to the wave operator in lower dimensions $n=1, 2, 3$. Recently, Ikeda and Ogawa \cite{Ikeda} established the lifespan estimate for Cauchy problem \eqref{1} for $n\geq 1$ in the form
 \begin{equation}\label{4}
\begin{aligned}
\exp\left(c\varepsilon^{-(p-1)}\right)\leq T(\varepsilon)\leq \exp\left(C\varepsilon^{-p}\right),\\
\end{aligned}
\end{equation}
obviously, there is a gap between the upper and lower bound of the lifespan estimate.

\begin{rem}
The test function method introduced in \cite{Zhang} is very useful for proving blow up result for both semilinear heat and damped wave equation,
also it is powerful for obtaining the sharp upper bound of lifespan for the corresponding semilinear problems with subcritical power. However, it doesn't work for
the critical power. In \cite{Ikeda}, in additional to the test function method, they also used a more delicate argument, by using the support
property of the test function to separate the time interval into small pieces $\thicksim 2^{j/2}$. Unfortunately, it seems we can't improve the upper bound of the lifespan estimate by the same argument.
\end{rem}

In this paper we are devoted to filling the gap between the upper bound and the lower bound of the lifespan estimate for problem \eqref{1} in higher dimensional spaces$(n\geq 4)$. As mentioned above, we can not use the positivity of the fundamental solution to the wave operator in higher dimensional case. However, we may use the idea of test function method for semilinear wave equations $u_{tt}-\Delta u=|u|^p$, which was introduced by Yordanov and Zhang \cite{Yordanov} and Zhou \cite{Zhou}.
The key ingredient now is to find an appropriate test function. Since the linear damped wave equation $u_{tt}-\Delta u+u_t=0$ has the so-called diffusion phenomenon and its solution behaves like that of the corresponding linear heat equation $v_t-\Delta v=0$ as $t\rightarrow \infty$, it will be useful to use the heat kernel, the fundamental solution of the heat equation, as the test function. Then we treat the equation in \eqref{1} as a heat equation, and hence we may get the expression of the solution in terms of the initial data and the source terms $-u_{tt}$ and $|u|^p$, by Duhamel's principle. After that we can establish an ordinary differential inequality for an appropriate functional to get the desired lifespan estimate. The essence is that we use the semigroup property of the heat kernel to deal with the source term $u_{tt}$.

\begin{thm}\label{thm1}
Let $f, g\in H^1(\mathbb{R}^n)\times L^2(\mathbb{R}^n)$and $supp (f, g)\subset \{x\in \mathbb{R}^n:\big| |x|\leq 1\}$. Assume furthermore that $f(x), g(x)\geq 0$ are non-trivial.
Let $u(t, x)$ solve the Cauchy problem \eqref{1} with $p=1+\frac2n$ on
$[0, T(\varepsilon))$. Then we have
\begin{equation}\label{5}
\begin{aligned}
T(\varepsilon)\leq \exp(C\varepsilon^{-\frac2n}),\\
\end{aligned}
\end{equation}
\end{thm}
where $C=C(n, f, g)$ is a positive constant independent of $\varepsilon$.

\begin{rem}
Our method also works for the lower dimensional cases $n=1, 2, 3$. The methods used in \cite{Li} for $n=1, 2$ and \cite{Nishihara} for $n=3$
include sophisticated analysis for the asymptotic behavior of some special functions, i.e., the modified Bessel function and hyperbolic cosine function.
\end{rem}

\section{Heat Kernel}

In this section we make a brief introduction to the heat kernel. In the $n$ dimensional case, the heat kernel in $\mathbb{R}^n$ is given by
\begin{equation}\label{5a}
\begin{aligned}
E(t, x)=\frac{1}{(4\pi t)^{\frac n2}}e^{-\frac{|x|^2}{4t}}.
\end{aligned}
\end{equation}
\begin{lem}[Evans \cite{Evans}]\label{lem1}
For any given time $t>0$, the heat kernel satisfies
\begin{equation}\label{5b}
\begin{aligned}
\int_{\mathbb{R}^n}E(t, x)dx=1,\\
\end{aligned}
\end{equation}
and
\begin{equation}\label{5c}
\left \{
\begin{aligned}
&E_t-\Delta E=0,~~~(t, x)\in (0, \infty)\times \mathbb{R}^n,\\
&E=\delta_{O},~~~on~ \{t=0\}\times \mathbb{R}^n,\\
\end{aligned} \right.
\end{equation}
where $\delta_{O}$ denotes the Dirac measure on $\mathbb{R}^n$, giving unit mass to the original point $O$.
\end{lem}

\begin{lem}[Semigroup property]\label{lem2}
For any given $t, s>0$ and $x\in \mathbb{R}^n$, the heat kernel $E(t, x)$ satisfies
\begin{equation}\label{5d}
\begin{aligned}
E(t, x)\ast E(s, x)=E(t+s, x),
\end{aligned}
\end{equation}
where $\ast$ denotes the convolution.
\end{lem}

Proof. By direction computation we know the Fourier transformation of $E(t, x)$ is given by
\begin{equation}\label{5e}
\begin{aligned}
\widehat{E}(t, \xi)=e^{-t|\xi|^2},
\end{aligned}
\end{equation}
and then we have
\begin{equation}\label{5f}
\begin{aligned}
\widehat{E}(t, \xi)\cdot \widehat{E}(s, \xi)=e^{-(t+s)|\xi|^2}=\widehat{E}(t+s, \xi),\\
\end{aligned}
\end{equation}
which implies \eqref{5d} in Lemma \ref{lem2}.

\section{Proof of the Main Theorem}

First the equation in Cauchy problem \eqref{1} can be regarded as the heat equation with source term
\begin{equation}\nonumber
\begin{aligned}
u_t-\Delta u=|u|^p-u_{tt},
\end{aligned}
\end{equation}
then by Duhamel's principle and using the fist initial data in \eqref{1} we have
\begin{equation}\label{6}
\begin{aligned}
u(t, x)=\varepsilon E(t)\ast f+\int_0^t E(t-\tau)\ast |u(\tau)|^pd\tau-\int_0^tE(t-\tau)\ast \partial_{\tau\tau}u(\tau)d\tau.
\end{aligned}
\end{equation}
By integration by parts and using the second initial data in \eqref{1} the last term can be rewritten as
\begin{equation}\label{7}
\begin{aligned}
&\int_0^tE(t-\tau)\ast \partial_{\tau\tau}u(\tau)d\tau\\
=&\int_{\mathbb{R}^n}\int_0^tE(t-\tau, x-y)\partial_{\tau\tau}u(\tau, y)d\tau dy\\
=&\int_{\mathbb{R}^n}E(0, x-y)u_t(y)dy-\int_{\mathbb{R}^n}E(t, x-y)u_t(0, y)dy\\
&-\int_{\mathbb{R}^n}\int_0^t\partial_{\tau}E(t-\tau, x-y)\partial_{\tau}u(\tau, y)d\tau dy\\
=&u_t(x)-\varepsilon E(t)\ast g-\int_{\mathbb{R}^n}\int_0^t\partial_{\tau}E(t-\tau, x-y)\partial_{\tau}u(\tau, y)d\tau dy,\\
\end{aligned}
\end{equation}
in which we used the second equality of \eqref{5c} in Lemma \ref{lem1}.
Plugging \eqref{7} into \eqref{6} we have
\begin{equation}\label{8}
\begin{aligned}
&u(t, x)+u_t(t, x)-\int_{\mathbb{R}^n}\int_0^t\partial_{\tau}E(t-\tau, x-y)\partial_{\tau}u(\tau, y)d\tau dy\\
=&\varepsilon E(t)\ast f+\varepsilon E(t)\ast g+\int_0^t E(t-\tau)\ast |u|^pd\tau.\\
\end{aligned}
\end{equation}
Multiplying the both sides of \eqref{8} with $\big(4\pi(t+1)\big)^{\frac n2}E(t+1)$ and then integrating over $\mathbb{R}^n$, we have
\begin{equation}\label{9}
\begin{aligned}
&\int_{\mathbb{R}^n}e^{-\frac{|x|^2}{4(t+1)}}u(t, x)dx+\int_{\mathbb{R}^n}e^{-\frac{|x|^2}{4(t+1)}}u_t(t, x)dx\\
&-\big(4\pi(t+1)\big)^{\frac n2}\int_{\mathbb{R}^n}\int_0^t\int_{\mathbb{R}^n}E(t+1, x)\partial_{\tau}E(t-\tau, x-y)\\
&~~\cdot\partial_{\tau}u(\tau, y)dyd\tau dx\\
=&\varepsilon \big(4\pi(t+1)\big)^{\frac n2}\int_{\mathbb{R}^n}E(2t+1, x)(f(x)+g(x))dx\\
&+\big(4\pi(t+1)\big)^{\frac n2}\int_0^{t}\int_{\mathbb{R}^n}E(2t+1-\tau, x)|u|^p(\tau, x)dxd\tau,\\
\end{aligned}
\end{equation}
in which we used the semigroup property of $E(t, x)$ given in Lemma \ref{lem2}.

Setting
\begin{equation}\label{10}
\begin{aligned}
G(t)&=\int_{\mathbb{R}^n}e^{-\frac{|x|^2}{4(t+1)}}u(t, x)dx,\\
F(t)&=\Big(\int_{\mathbb{R}^n}e^{-\frac{|x|^2}{4(t+1)}}|u(t, x)|^pdx\Big)^{\frac1p}(t+1)^{\frac{n(p-1)}{2p}}\\
&=\Big(\int_{\mathbb{R}^n}e^{-\frac{|x|^2}{4(t+1)}}|u(t, x)|^pdx\Big)^{\frac1p}(t+1)^{\frac{n}{n+2}}\\
\end{aligned}
\end{equation}
and
\begin{equation}\label{10a}
\begin{aligned}
A(t)=&\int_{\mathbb{R}^n}e^{-\frac{|x|^2}{4(t+1)}}u_t(t, x)dx,\\
B(t)=&-\big(4\pi(t+1)\big)^{\frac n2}\int_{\mathbb{R}^n}\int_0^t\int_{\mathbb{R}^n}E(t+1, x)\partial_{\tau}E(t-\tau, x-y)\\
&~~~\cdot \partial_{\tau}u(\tau, y)dyd\tau dx,\\
D(t)=&\big(4\pi(t+1)\big)^{\frac n2}\int_0^{t}\int_{\mathbb{R}^n}E(2t+1-\tau, x)|u|^p(\tau, x)dxd\tau,\\
\end{aligned}
\end{equation}
by H\"{o}lder inequality we have
\begin{equation}\label{11}
\begin{aligned}
G(t)&=\int_{\mathbb{R}^n}e^{-\frac{|x|^2}{4(t+1)}}u(t, x)dx\\
&\leq \Big(\int_{\mathbb{R}^n}e^{-\frac{|x|^2}{4(t+1)}}|u(t, x)|^pdx\Big)^{\frac1p}\Big(\int_{\mathbb{R}^n}e^{-\frac{|x|^2}{4(t+1)}}dx\Big)^{\frac{1}{p'}}\\
&\leq C\Big(\int_{\mathbb{R}^n}e^{-\frac{|x|^2}{4(t+1)}}|u(t, x)|^pdx\Big)^{\frac1p}(t+1)^{\frac{n(p-1)}{2p}}\Big(\int_{\mathbb{R}^n}e^{-|y|^2}dy\Big)^{\frac{1}{p'}}\\
&\leq CF(t).\\
\end{aligned}
\end{equation}
Here and hereafter, $C$ denotes a positive constant which is independent of $\varepsilon$ and may change from line to line. As to $A(t)$, a direct
computation leads to
\begin{equation}\label{12}
\begin{aligned}
A(t)&=\int_{\mathbb{R}^n}e^{-\frac{|x|^2}{4(t+1)}}u_t(t, x)dx\\
&=\frac{d}{dt}\int_{\mathbb{R}^n}e^{-\frac{|x|^2}{4(t+1)}}u(t, x)dx-\int_{\mathbb{R}^n}e^{-\frac{|x|^2}{4(t+1)}}u(t, x)\frac{|x|^2}{4(t+1)^2}dx\\
\end{aligned}
\end{equation}
with
\begin{equation}\label{13}
\begin{aligned}
&\int_{\mathbb{R}^n}e^{-\frac{|x|^2}{4(t+1)}}u(t, x)\frac{|x|^2}{4(t+1)^2}dx\\
\leq& \Big(\int_{\mathbb{R}^n}e^{-\frac{|x|^2}{4(t+1)}}|u(t, x)|^pdx\Big)^{\frac1p}\Big(\int_{\mathbb{R}^n}e^{-\frac{|x|^2}{4(t+1)}}\big(\frac{|x|^2}{4(t+1)^2}\big)^{p'}dx\Big)^{\frac{1}{p'}}\\
\leq& C(t+1)^{-1}F(t),\\
\end{aligned}
\end{equation}
hence
\begin{equation}\label{14}
\begin{aligned}
A(t)\leq G'(t)+C(t+1)^{-1}F(t).\\
\end{aligned}
\end{equation}

It is easy to see that $2(t+1)\geq 2t+1-\tau\geq \tau+1$ for $0\leq \tau\leq t$, then the nonlinear term $D(t)$ can be estimated by
\begin{equation}\label{15}
\begin{aligned}
D(t)&=\big(4\pi(t+1)\big)^{\frac n2}\int_0^{t}\int_{\mathbb{R}^n}E(2t+1-\tau, x)|u|^p(\tau, x)dxd\tau\\
&=\int_0^{t}\frac{\big(4\pi(t+1)\big)^{\frac n2}}{\big(4\pi(2t+1-\tau)\big)^{\frac n2}}\int_{\mathbb{R}^n}e^{-\frac{|x|^2}{4(2t+1-\tau)}}|u|^pdxd\tau\\
&\geq 2^{-\frac n2}\int_0^{t}\int_{\mathbb{R}^n}e^{-\frac{|x|^2}{4(\tau+1)}}|u|^pdxd\tau\\
&\geq C\int_0^{t}\frac{F^p(\tau)}{\tau+1}d\tau,\\
\end{aligned}
\end{equation}
in which we used the definition of $F(\tau)$ given by \eqref{10}. Therefore, the combination of \eqref{9}-\eqref{10a} and \eqref{14} yields
\begin{equation}\label{16}
\begin{aligned}
&G'(t)+G(t)+\frac{CF(t)}{t+1}+B(t)\\
\geq &\varepsilon \big(4\pi(t+1)\big)^{\frac n2}\int_{\mathbb{R}^n}E(2t+1, x)\big(f(x)+g(x)\big))dx+D(t)\\
=&\varepsilon \Big(\frac{t+1}{2t+1}\Big)^{\frac n2}\int_{\mathbb{R}^n}e^{-\frac{|x|^2}{4(2t+1)}}\big(f(x)+g(x)\big)dx+D(t).\\
\end{aligned}
\end{equation}

We now estimate the term $B(t)$. By integration by parts and using Lemma \ref{lem2} we have
\begin{equation}\label{17}
\begin{aligned}
&B(t)\\
=&-\big(4\pi(t+1)\big)^{\frac n2}\int_{\mathbb{R}^n}\int_0^t\int_{\mathbb{R}^n}E(t+1, x)\partial_{\tau}E(t-\tau, x-y)\\
&~~\cdot\partial_{\tau}u(\tau, y)dyd\tau dx\\
=&-\big(4\pi(t+1)\big)^{\frac n2}\int_0^t\int_{\mathbb{R}^n}\partial_{\tau}\big(\int_{\mathbb{R}^n}E(t+1, x)E(t-\tau, x-y)dx\big)\\
&~~\cdot \partial_{\tau}u(\tau, y)dyd\tau\\
=&-\big(4\pi(t+1)\big)^{\frac n2}\int_0^t\int_{\mathbb{R}^n}\partial_{\tau}E(2t+1-\tau, x)\partial_{\tau}u(\tau, x)dxd\tau\\
=&-\big(4\pi(t+1)\big)^{\frac n2}\int_0^t\partial_{\tau}\big(\int_{\mathbb{R}^n}E_{\tau}(2t+1-\tau, x)u(\tau, x)dx\big)d\tau\\
&+\big(4\pi(t+1)\big)^{\frac n2}\int_0^t\int_{\mathbb{R}^n}\partial_{\tau\tau}E(2t+1-\tau, x)u(\tau, x)dxd\tau\\
\triangleq& B_1(t)+B_2(t).\\
\end{aligned}
\end{equation}
There are four terms in $B_1(t)$:
\begin{equation}\label{18}
\begin{aligned}
B_1(t)=&-\big(4\pi(t+1)\big)^{\frac n2}\Big[\int_{\mathbb{R}^n}\frac n2\big(4\pi(t+1)\big)^{-\frac n2-1}e^{-\frac{|x|^2}{4(t+1)}}u(t, x)dx\\
&-\int_{\mathbb{R}^n}\big(4\pi(t+1)\big)^{-\frac n2}e^{-\frac{|x|^2}{4(t+1)}}u(t, x)\frac{|x|^2}{4(t+1)^2}dx\Big]\\
&+\big(4\pi(t+1)\big)^{\frac n2}\Big[\int_{\mathbb{R}^n}\frac n2\big(4\pi(2t+1)\big)^{-\frac n2-1}e^{-\frac{|x|^2}{4(2t+1)}}u(0, x)dx\\
&-\int_{\mathbb{R}^n}\big(4\pi(2t+1)\big)^{-\frac n2}e^{-\frac{|x|^2}{4(2t+1)}}u(0, x)\frac{|x|^2}{4(2t+1)^2}dx\Big]\\
\triangleq& B_{11}(t)+B_{12}(t)+B_{13}(t)+B_{14}(t).\\
\end{aligned}
\end{equation}
By \eqref{11} and \eqref{13} it is easy to get
\begin{equation}\label{19}
\begin{aligned}
B_{11}(t)&=-\big(4\pi(t+1)\big)^{\frac n2}\int_{\mathbb{R}^n}\frac n2\big(4\pi(t+1)\big)^{-\frac n2-1}e^{-\frac{|x|^2}{4(t+1)}}u(t, x)dx\\
&\leq C(t+1)^{-1}|G(t)|\\
&\leq C(t+1)^{-1}F(t),\\
B_{12}(t)&=\int_{\mathbb{R}^n}e^{-\frac{|x|^2}{4(t+1)}}u(t, x)\frac{|x|^2}{4(t+1)^2}dx\\
&\leq C(t+1)^{-1}F(t).\\
\end{aligned}
\end{equation}
The other two terms $B_{13}$ and $B_{14}$ are related to the initial data:
\begin{equation}\label{20}
\begin{aligned}
B_{13}(t)&=\varepsilon \frac n2\frac{1}{4\pi(2t+1)}\Big(\frac{t+1}{2t+1}\Big)^{\frac n2}\int_{\mathbb{R}^n}e^{-\frac{|x|^2}{4(2t+1)}}f(x)dx,\\
B_{14}(t)&=-\varepsilon \Big(\frac{t+1}{2t+1}\Big)^{\frac n2}\int_{\mathbb{R}^n}e^{-\frac{|x|^2}{4(2t+1)}}f(x)\frac{|x|^2}{4(2t+1)^2}dx.\\
\end{aligned}
\end{equation}
It follows then from \eqref{19} and \eqref{20} that
\begin{equation}\label{20a}
\begin{aligned}
B_{1}(t)&\leq C(t+1)^{-1}F(t)+\varepsilon \frac n2\frac{1}{4\pi(2t+1)}\Big(\frac{t+1}{2t+1}\Big)^{\frac n2}\int_{\mathbb{R}^n}e^{-\frac{|x|^2}{4(2t+1)}}f(x)dx\\
&-\varepsilon \Big(\frac{t+1}{2t+1}\Big)^{\frac n2}\int_{\mathbb{R}^n}e^{-\frac{|x|^2}{4(2t+1)}}f(x)\frac{|x|^2}{4(2t+1)^2}dx.\\
\end{aligned}
\end{equation}
There are also four terms in $B_2(t)$:
\begin{equation}\label{21}
\begin{aligned}
&B_2(t)\\
=&\big(4\pi(t+1)\big)^{\frac n2}\int_0^t\int_{\mathbb{R}^n}\partial_{\tau\tau}E(2t+1-\tau, x)u(\tau, x)dxd\tau\\
=&\big(4\pi(t+1)\big)^{\frac n2}\Big[\int_0^t\int_{\mathbb{R}^n}\frac n2(\frac n2+1)\big(4\pi(2t+1-\tau)\big)^{-\frac n2-2}\\
 &\cdot e^{-\frac{|x|^2}{4(2t+1-\tau)}}
u(\tau, x)dxd\tau\\
-&\int_0^t\int_{\mathbb{R}^n}n\big(4\pi(2t+1-\tau)\big)^{-\frac n2-1}e^{-\frac{|x|^2}{4(2t+1-\tau)}}\frac{|x|^2}{4(2t+1-\tau)^2}
u(\tau, x)dxd\tau\\
+&\int_0^t\int_{\mathbb{R}^n}\big(4\pi(2t+1-\tau)\big)^{-\frac n2}e^{-\frac{|x|^2}{4(2t+1-\tau)}}\frac{|x|^4}{16(2t+1-\tau)^4}
u(\tau, x)dxd\tau\\
-&\int_0^t\int_{\mathbb{R}^n}\big(4\pi(2t+1-\tau)\big)^{-\frac n2}e^{-\frac{|x|^2}{4(2t+1-\tau)}}\frac{|x|^2}{2(2t+1-\tau)^3}
u(\tau, x)dxd\tau\Big]\\
\triangleq &B_{21}(t)+B_{22}(t)+B_{23}(t)+B_{24}(t).\\
\end{aligned}
\end{equation}
Since we have $2t+1-\tau\geq t+1$ for $0\leq \tau\leq t$, by H\"{o}lder inequality
and Yong inequality we get
\begin{equation}\label{22}
\begin{aligned}
B_{21}(t)\leq &C(t+1)^{\frac n2-2}\int_0^t\int_{\mathbb{R}^n}E(2t+1-\tau, x)u(\tau, x)dxd\tau\\
\leq &C(t+1)^{\frac n2-2}\Big(\int_0^t\int_{\mathbb{R}^n}E(2t+1-\tau, x)|u|^p(\tau, x)dxd\tau\Big)^{\frac1p}\\
&\cdot \Big(\int_0^t\int_{\mathbb{R}^n}E(2t+1-\tau, x)dxd\tau\Big)^{\frac1{p'}}\\
\leq &\frac18\big(4\pi(t+1)\big)^{\frac n2}\int_0^t\int_{\mathbb{R}^n}E(2t+1-\tau, x)|u|^p(\tau, x)dxd\tau\\
&+C(t+1)^{-\frac n2-1}\\
=&\frac18 D(t)+C(t+1)^{-\frac n2-1}\\
\end{aligned}
\end{equation}
and
\begin{equation}\label{23}
\begin{aligned}
B_{22}(t)\leq& C(t+1)^{\frac n2-1}\int_0^t\int_{\mathbb{R}^n}E(2t+1-\tau, x)\frac{|x|^2}{4(2t+1-\tau)^2}u(\tau, x)dxd\tau\\
\leq &C(t+1)^{\frac n2-1}\Big(\int_0^t\int_{\mathbb{R}^n}E(2t+1-\tau, x)|u|^p(\tau, x)dxd\tau\Big)^{\frac1p}\\
&\cdot \Big(\int_0^t\int_{\mathbb{R}^n}E(2t+1-\tau, x)\Big(\frac{|x|^2}{4(2t+1-\tau)^2}\Big)^{p'}dxd\tau\Big)^{\frac1{p'}}\\
\leq &C(t+1)^{\frac n2-1}\Big(\int_0^t\int_{\mathbb{R}^n}E(2t+1-\tau, x)|u|^p(\tau, x)dxd\tau\Big)^{\frac1p}\\
&\cdot \Big(\int_0^t\frac{1}{(2t+1-\tau)^{p'}}\int_{\mathbb{R}^n}e^{-|y|^2}|y|^{2p'}dyd\tau\Big)^{\frac1{p'}}\\
\leq&\frac18 D(t)+C(t+1)^{-\frac n2-1},\\
\end{aligned}
\end{equation}
in which we used the variable transformation $x=2\sqrt{2t+1-\tau}y$ in the third inequality of \eqref{23}. By a similar way to that of $B_{22}$, we have
\begin{equation}\label{24}
\begin{aligned}
&B_{23}(t)\leq \frac18 D(t)+C(t+1)^{-\frac n2-1},\\
&B_{24}(t)\leq \frac18 D(t)+C(t+1)^{-\frac n2-1}.\\
\end{aligned}
\end{equation}
By combining \eqref{22}-\eqref{24} we get
\begin{equation}\label{25}
\begin{aligned}
&B_{2}(t)\leq \frac12 D(t)+C(t+1)^{-\frac n2-1},\\
\end{aligned}
\end{equation}
which together with \eqref{20a} gives
\begin{equation}\label{26}
\begin{aligned}
B(t)=&B_1(t)+B_2(t)\\
\leq &C(t+1)^{-1}F(t)+\varepsilon \frac n2\frac{1}{4\pi(2t+1)}\Big(\frac{t+1}{2t+1}\Big)^{\frac n2}\int_{\mathbb{R}^n}e^{-\frac{|x|^2}{4(2t+1)}}f(x)dx\\
&-\varepsilon \Big(\frac{t+1}{2t+1}\Big)^{\frac n2}\int_{\mathbb{R}^n}e^{-\frac{|x|^2}{4(2t+1)}}f(x)\frac{|x|^2}{4(2t+1)^2}dx\\
&+\frac12 D(t)+C_0(t+1)^{-\frac n2-1},\\
\end{aligned}
\end{equation}
where $C_0>0$ is a constant depending only on $(\pi, n)$. By \eqref{15}-\eqref{16} and \eqref{26} we get
\begin{equation}\label{27}
\begin{aligned}
&G'(t)+G(t)+\frac{CF(t)}{t+1}\\
\geq& \frac12 D(t)+\big(1-\frac n2\frac1{4\pi(2t+1)}\big)\varepsilon \Big(\frac{t+1}{2t+1}\Big)^{\frac n2}\int_{\mathbb{R}^n}e^{-\frac{|x|^2}{4(2t+1)}}f(x)dx\\
&+\varepsilon \Big(\frac{t+1}{2t+1}\Big)^{\frac n2}\int_{\mathbb{R}^n}e^{-\frac{|x|^2}{4(2t+1)}}f(x)\frac{|x|^2}{4(2t+1)^2}dx\\
&+\varepsilon \Big(\frac{t+1}{2t+1}\Big)^{\frac n2}\int_{\mathbb{R}^n}e^{-\frac{|x|^2}{4(2t+1)}}g(x)dx-C_0(t+1)^{-\frac n2-1}\\
\geq& C_1\int_0^{t}\frac{F^p(\tau)}{\tau+1}d\tau+\big(1-\frac n2\frac1{4\pi(2t+1)}\big)\varepsilon \Big(\frac{t+1}{2t+1}\Big)^{\frac n2}\int_{\mathbb{R}^n}e^{-\frac{|x|^2}{4(2t+1)}}f(x)dx\\
&+\varepsilon \Big(\frac{t+1}{2t+1}\Big)^{\frac n2}\int_{\mathbb{R}^n}e^{-\frac{|x|^2}{4(2t+1)}}f(x)\frac{|x|^2}{4(2t+1)^2}dx\\
&+\varepsilon \Big(\frac{t+1}{2t+1}\Big)^{\frac n2}\int_{\mathbb{R}^n}e^{-\frac{|x|^2}{4(2t+1)}}g(x)dx-C_0(t+1)^{-\frac n2-1},\\
\end{aligned}
\end{equation}
where $C_1>0$ is a constant depending only on $n$ as it comes from \eqref{15}.

Setting
\begin{equation}\nonumber
\begin{aligned}
H(t)=&\Big(1-\frac n2\frac1{4\pi(2t+1)}\Big)\Big(\frac{t+1}{2t+1}\Big)^{\frac n2}\int_{\mathbb{R}^n}e^{-\frac{|x|^2}{4(2t+1)}}f(x)dx\\
&+\Big(\frac{t+1}{2t+1}\Big)^{\frac n2}\int_{\mathbb{R}^n}e^{-\frac{|x|^2}{4(2t+1)}}f(x)\frac{|x|^2}{4(2t+1)^2}dx\\
&+\Big(\frac{t+1}{2t+1}\Big)^{\frac n2}\int_{\mathbb{R}^n}e^{-\frac{|x|^2}{4(2t+1)}}g(x)dx\\
\end{aligned}
\end{equation}
and
\begin{equation}\nonumber
\begin{aligned}
I(t)=C_0(t+1)^{-\frac n2-1}.\\
\end{aligned}
\end{equation}
Multiplying \eqref{27} with $e^t$ and then integrating it over $[0, t]$, we have
\begin{equation}\label{27a}
\begin{aligned}
e^tF(t)+\int_{0}^{t}\frac{e^{\tau}F(\tau)}{\tau+1}d\tau\geq &e^tG(t)+\int_{0}^{t}\frac{e^{\tau}F(\tau)}{\tau+1}d\tau\\
\geq & G(0)+C_1\int_{0}^{t}\frac{(e^t-e^{\tau})F^p(\tau)}{\tau+1}d\tau\\
&+\int_0^te^{\tau}\big(\varepsilon H(\tau)-I(\tau)\big)d\tau.\\
\end{aligned}
\end{equation}
We claim that
\begin{equation}\label{27b}
\begin{aligned}
J(t):=\int_0^te^{\tau}\big(\varepsilon H(\tau)-I(\tau)\big)d\tau\geq C_2\e e^t
\end{aligned}
\end{equation}
for some constant $C_2>0$ depending only on $(f, g, n)$ and $t$ large enough. Actually, since we assume that the initial data $f(x), g(x)\in L^2(\R^n)$ and have compact support, we have
$f(x), g(x)\in L^1(\R^n)$. So by the dominated convergence theorem we have
\begin{equation}\label{27c}
\begin{aligned}
\lim_{t\rightarrow \infty}H(t)=2^{-\frac n2}\int_{\R^n}\big(f(x)+g(x)\big)dx=H_0>0,
\end{aligned}
\end{equation}
and then there exists $t_1=t_1(n, f, g)$ independent of $\varepsilon$, such that
\begin{equation}\label{27d}
\begin{aligned}
H(t)>\frac{H_0}{2}>0,~~~for~t>t_1.
\end{aligned}
\end{equation}
Let $t_2=\big(\frac{H_0}{4C_0}\big)^{-\frac{2}{n+2}}\e^{-\frac{2}{n+2}}-1$. We get
\begin{equation}\label{27e}
\begin{aligned}
I(t)\leq \frac{H_0\e}{4},~~~for~t\geq t_2.\\
\end{aligned}
\end{equation}
Setting $t_3=\max\{t_1, t_2\}$, since $H(t)$ is continuous and has a limit as $t\rightarrow \infty$, $H(t)$ is bounded over $[0, \infty)$, and then we may set $M_{H}=\max_{t\in [0, \infty)}|H(t)|<\infty$. On the other hand, it is easy to get $M_{I}=\max_{t\in [0, \infty)}I(t)=C_0>0$. Hence we have
\begin{equation}\label{27f}
\begin{aligned}
J(t)&:=\int_0^te^{\tau}\big(\e H(\tau)-I(\tau)\big)d\tau\\
&=\int_0^{t_3}e^{\tau}\big(\e H(\tau)-I(\tau)\big)d\tau+\int_{t_3}^te^{\tau}\big(\e H(\tau)-I(\tau)\big)d\tau\\
&\geq \int_0^{t_3}e^{\tau}\big(-\e M_H-C_0\big)d\tau+\frac{H_0\e}{4}(e^t-e^{t_3})\\
&\geq \frac{H_0\e}{8}e^t+\frac{H_0\e}{8}e^t-\big(\frac{H_0\e}{4}+M_H \e+C_0\big)e^{t_3}+M_H\e+C_0.\\
\end{aligned}
\end{equation}
Hence, if $t\geq t_4=\ln\big(2+\frac{8M_H}{H_0}+\frac{8C_0}{H_0}\e^{-1}\big)+t_3$, then we get
\begin{equation}\label{27g}
\begin{aligned}
J(t)\geq \frac{H_0\e}{8}e^t+M_H\e+C_0\geq C_2\e e^t,
\end{aligned}
\end{equation}
which proves \eqref{27b}. With this claim in mind, it is easy to get from \eqref{27a} that
\begin{equation}\label{27h}
\begin{aligned}
e^tF(t)+\int_{0}^{t}\frac{e^{\tau}F(\tau)}{\tau+1}d\tau
\geq C_2\e e^t+C_1\int_{0}^{t}\frac{(e^t-e^{\tau})F^p(\tau)}{\tau+1}d\tau~~~for~t\geq t_4,\\
\end{aligned}
\end{equation}
which can be rewritten as
\begin{equation}\label{33}
\begin{aligned}
\Big((1+t)\int_{0}^{t}\frac{e^{\tau}F(\tau)}{\tau+1}d\tau\Big)'\geq C_2\varepsilon e^t+C_1\int_{0}^{t}\frac{(e^t-e^{\tau})F^p(\tau)}{\tau+1}d\tau~~for~t\geq t_4.
\end{aligned}
\end{equation}
Setting $t_0=t_4+\ln2$ and
\[
\begin{aligned}
\alpha(t)&=(1+t)\int_{0}^{t}\frac{e^{\tau}F(\tau)}{\tau+1}d\tau,\\
\beta(t)&=\int_{t_0}^t\int_{0}^s\frac{(e^s-e^{\tau})F^p(\tau)}{1+\tau}
d\tau ds.\\
\end{aligned}
\]
Hence integrating \eqref{33} over $[t_4, t]$, we obtain
\begin{equation}\label{34}
\begin{aligned}
\alpha(t)-\alpha(t_4)\geq C_2\varepsilon (e^t-e^{t_4})+C_1\left(\beta(t)-\beta(t_4)\right),~~~~~t\geq t_4.\\
\end{aligned}
\end{equation}
Noting that $\alpha(t_4)\geq 0$ and $\beta(t_4)\leq 0$, we then get from the above inequality that
\begin{equation}\label{34a}
\begin{aligned}
\alpha(t)&\geq \frac12C_2\e e^t+\frac12C_2\e e^t-C_2\e e^{t_4}+C_1\beta(t)\\
&\geq \frac12C_2\e e^t+C_1\beta(t)\quad\mbox{for}\ t\ge t_0.
\end{aligned}
\end{equation}
Also from the definition of $\beta(t)$ we have
\begin{equation}\nonumber
\begin{aligned}
\beta'(t)=\int_{0}^{t}\frac{(e^t-e^{\tau})F^p(\tau)}{\tau+1}d\tau,\\
\end{aligned}
\end{equation}
then
\begin{equation}\label{35}
\begin{aligned}
\beta''(t)-\beta'(t)=\int_{0}^{t}\frac{e^{\tau}F^p(\tau)}{\tau+1}d\tau.\\
\end{aligned}
\end{equation}

Before going on, we have another claim that
\begin{equation}\label{35a}
\begin{aligned}
\int_0^t\frac{e^{\tau}}{1+\tau}d\tau\leq C_3\frac{e^t}{1+t},~~~for~t\geq 0,
\end{aligned}
\end{equation}
where $C_3$ denotes a positive constant independent of $\e$. For this claim, first it is easy to get that there exists a positive constant $\widetilde{C}_3$ independent of $\e$ such that
\[
\int_{\frac t2}^t\frac{e^{\tau}}{1+\tau}d\tau\leq \widetilde{C}_3\frac{e^t-e^{\frac t2}}{1+t}\leq \widetilde{C}_3\frac{e^t}{1+t},~~~for~t\geq 0.
\]
Since the integrand $K(t)=\frac{e^t}{1+t}$ is an increasing function of $t$, then we have
\[
\int_0^{\frac t2}\frac{e^{\tau}}{1+\tau}d\tau\leq \int_{\frac t2}^t\frac{e^{\tau}}{1+\tau}d\tau\leq \widetilde{C}_3\frac{e^t}{1+t},~~~for~t\geq 0,
\]
hence we prove claim \eqref{35a} by setting $C_3=2\widetilde{C}_3$. By H\"{o}lder inequality and \eqref{35a} we obtain that
\begin{equation}\label{36}
\begin{aligned}
\alpha(t)&\leq (1+t)\Big(\int_{0}^{t}\frac{e^{\tau}F^p(\tau)}{\tau+1}d\tau\Big)^{\frac1p}\Big(\int_{0}^t\frac{e^{\tau}}{1+\tau}d\tau\Big)^{\frac1{p'}}\\
&\leq C_4e^{\frac t{p'}}(1+t)^{\frac1p}\Big(\int_{0}^{t}\frac{e^{\tau}F^p(\tau)}{\tau+1}d\tau\Big)^{\frac1p},\\
\end{aligned}
\end{equation}
where $C_4:=(C_3)^{1/p'}$. Thus, combining \eqref{34a}-\eqref{35} and \eqref{36} we have
\begin{equation}\label{37}
\begin{aligned}
\beta''(t)-\beta'(t)&\geq \frac{C_5\alpha^p(t)e^{-(p-1)t}}{1+t}\\
&\geq \frac{C_6\big(\varepsilon e^t+\beta(t)\big)^pe^{-(p-1)t}}{1+t},~~~t\geq t_0.\\
\end{aligned}
\end{equation}
where $C_5:=(C_4)^{-p}$ and
\[
C_6:=C_5\left(\min\{\frac12 C_2, C_1\}\right)^p.\\
\]
Setting $\beta(t)=e^t\gamma(t)$ and plugging it into \eqref{37}, we get
\begin{equation}\label{38}
\begin{aligned}
\gamma''(t)+\gamma'(t)\geq \frac{C_6\big(\varepsilon+\gamma(t)\big)^p}{1+t},~~~t\geq t_0.\\
\end{aligned}
\end{equation}
Then, setting
\begin{equation}\label{50}
\begin{aligned}
\widetilde{\gamma}(t)=\varepsilon+\gamma(t)=\varepsilon+e^{-t}\beta(t)=\varepsilon+e^{-t}\int_{t_0}^t\int_{0}^s\frac{(e^s-e^{\tau})F^p(\tau)}{1+\tau}
d\tau ds\\
\end{aligned}
\end{equation}
with
\begin{equation}\label{51}
\begin{aligned}
\widetilde{\gamma}(t_0)=\varepsilon>0\\
\end{aligned}
\end{equation}
and
\begin{equation}\label{51a}
\begin{aligned}
\widetilde{\gamma}'(t_0)>0,\\
\end{aligned}
\end{equation}
it follows from \eqref{38} that
\begin{equation}\label{39}
\begin{aligned}
\widetilde{\gamma}''(t)+\widetilde{\gamma}'(t)\geq \frac{C_6\widetilde{\gamma}^p(t)}{1+t},~~~~~t\geq t_0.\\
\end{aligned}
\end{equation}

We have
\begin{lem}[Theorem 3.1 of Li and Zhou \cite{Li}]\label{lem3}
Suppose that $I(t)$ satisfies
\begin{equation}\label{52}
\begin{aligned}
I''(t)+I'(t)\geq C_0\frac{I^{1+\alpha}(t)}{(1+t)^{\beta}}~~(C_0>0, constant)
\end{aligned}
\end{equation}
and
\begin{equation}\label{53}
\begin{aligned}
I(0)>0,~~~I'(0)\geq 0,
\end{aligned}
\end{equation}
where $\alpha>0$. Then, when $0\leq \beta\leq 1$, $I=I(t)$ must blow up in a finite time. Moreover, if $I(0)=\varepsilon$, where $\varepsilon>0$ is a small parameter, then the lifespan $T(\varepsilon)$ of $I=I(t)$ has the following upper bound:
\begin{equation}\label{55}
T(\varepsilon)\leq\left \{
\begin{aligned}
&\exp(a\varepsilon^{-\alpha}),~~~~~\beta=1,\\
&b\varepsilon^{-\frac{\alpha}{1-\beta}},~~~~~0\leq \beta<1,\\
\end{aligned} \right.
\end{equation}
where $a$ and $b$ are positive constants independent of $\varepsilon$.
\end{lem}

Applying Lemma \ref{lem3} to $I(t)=\widetilde{\gamma}(t)$
with $p=1+\alpha=1+\frac2n$ and taking $\beta=1$ and the initial time $t=t_0$,
we then get the desired sharp upper bound of the lifespan given in Theorem \ref{thm1}.

\section*{Acknowledgement}
The first author is partially supported by Zhejiang Province
Science Foundation(LY18A010008), NSFC(11501273, 11726612, 11771359,
11771194), Chinese Postdoctoral Science Foundation(2017M620128), the Scientific Research Foundation of the First-Class Discipline of Zhejiang Province
(B)(201601). The second author is supported by Key Laboratory of Mathematics for Nonlinear
Sciences(Fudan University), Ministry of Education of China, Shanghai Key Laboratory for Contemporary
Applied Mathematics, School of Mathematical Sciences, Fudan University, NSFC (11726611, 11421061), 973 program
(2013CB834100) and 111 project.

The authors want to express their sincere thanks to Prof. T.T. Li for his helpful suggestion, they are also grateful to Prof. H. Takamura and Prof. M.Ikeda for pointing out the missing integral constant in \eqref{27a}.

\end{document}